\documentclass[letterpaper, 10 pt, conference]{ieeeconf}  

\IEEEoverridecommandlockouts                              

\usepackage{times} 
\usepackage{amsmath} 
\usepackage{amssymb}  

\usepackage{frenchineq}
\usepackage{amsfonts}
\usepackage[colorinlistoftodos,bordercolor=orange,backgroundcolor=orange!20,linecolor=orange,textsize=scriptsize]{todonotes}
\usepackage{hyperref}
\hypersetup{colorlinks=true,urlcolor=blue,linkcolor=blue,citecolor=red}
\usepackage{tikz}
\usetikzlibrary{shapes,arrows,decorations.markings,positioning}
\usepackage{hypergraph}

\newtheorem{theorem}{Theorem}

\newtheorem{proposition}[theorem]{Proposition}
\newtheorem{corollary}[theorem]{Corollary}
\newtheorem{definition}[theorem]{Definition}

\renewcommand{\theenumi}{\roman{enumi}}

\newcommand{\MIN}{{\small \textsf{MIN}}}
\newcommand{\MAX}{{\small \textsf{MAX}}}
\newcommand{\R}{\mathbb{R}}
\newcommand{\N}{\mathbb{N}}
\newcommand{\state}{S}
\newcommand{\Gcal}{\mathcal{G}}
\newcommand{\Hcal}{\mathcal{H}}
\newcommand{\Scal}{\mathcal{S}}

\newcommand{\unit}{e}
\newcommand{\tail}{\mathbf{t}}
\newcommand{\head}{\mathbf{h}}
\newcommand{\cpt}[1]{{\state \setminus #1}}
\newcommand{\Hnorm}[1]{\| \ifx\\#1\\ \cdot \else #1 \fi \|_{\mathsf{H}}}

\DeclareMathOperator{\size}{size}
\DeclareMathOperator{\reach}{reach}

\title{\LARGE \bf
Hypergraph conditions for the solvability of the ergodic equation for zero-sum games
}

\author{Marianne Akian$^{* 1}$, St\'ephane Gaubert$^{* 2}$ and Antoine Hochart$^{* \dag 3}$
\thanks{$^{*}$M.~Akian, S.~Gaubert and A.~Hochart are with INRIA Saclay-Ile-de-France and CMAP Ecole polytechnique, Route de Saclay, 91128 Palaiseau Cedex, France}
\thanks{$^\dag${A.~Hochart is supported by a PhD fellowship of Fondation Math\'ematique Jacques Hadamard (FMJH)}}%
\thanks{$^1${\tt\small marianne.akian@inria.fr}}
\thanks{$^2${\tt\small stephane.gaubert@inria.fr}}%
\thanks{$^3${\tt\small antoine.hochart@cmap.polytechnique.fr}}%
}

\begin{document}

\maketitle
\thispagestyle{empty}
\pagestyle{empty}

\begin{abstract}
The ergodic equation is a basic tool in the study of mean-payoff stochastic games.
Its solvability entails that the mean payoff is independent of the initial state.
Moreover, optimal stationary strategies are readily obtained from its solution.
In this paper, we give a general sufficient condition for the solvability of the ergodic equation, for a game with finite state space but arbitrary action spaces.
This condition involves a pair of directed hypergraphs depending only on the ``growth at infinity'' of the Shapley operator of the game.
This refines a recent result of the authors which only applied to games with bounded payments, as well as earlier nonlinear fixed point results for order preserving maps, involving graph conditions.
\end{abstract}

\begin{keywords}
Zero-sum games, stochastic control, ergodic control, risk-sensitive control, nonlinear consensus, computational methods, directed hypergraphs.
\end{keywords}

\section{Introduction}
\label{sec:Intro}

\paragraph{Motivation}
A general issue in the study of stochastic control and zero-sum game problems is to give conditions which guarantee that the mean payoff per time unit is independent of the choice of the initial state.
A basic tool to address this issue is the ergodic eigenproblem.
The latter is a nonlinear equation, involving a function or vector called {\em bias} or {\em potential}, and a scalar called {\em ergodic constant}.
When the ergodic equation has a solution, the ergodic constant gives the mean payoff for every initial state. 
The ergodic equation is often called the {\em average case optimality equation} in stochastic control, see~\cite{Whi86,HLL99}.
It may be thought of as a nonlinear extension of the {\em Poisson equation} arising in potential theory.

Conditions which guarantee the solvability of the ergodic equation are generally referred to as {\em ergodicity conditions}.
Such conditions typically involve the transition probabilities. 
One of the oldest conditions of this kind was obtained by Bather~\cite{Bat73}, who showed that the ergodic equation has a solution for the class of Markov decision processes that are {\em communicating}, meaning
that for every pair of states $i$ to $j$, there is a strategy and a horizon $k$ under which the probability to reach $j$ from $i$ is positive.
An interest of Bather's condition is its structural nature:
as it only depends on the transition probabilities, it is insensitive to a perturbation of the payments.
However, a limitation of Bather's communication condition is that it is tied to the one-player case:
the communication condition is almost never satisfied for two-player problems (including zero-sum deterministic games).

The question of finding generalized communication or ergodicity conditions has emerged in a number of works dealing with control and games~\cite{Kol92,KM97,CCHH09,CCHH10,AGH15}, but also discrete event systems~\cite{Ols91,YZ04}, mathematical economy~\cite{Osh83}, and Perron-Frobenius theory~\cite{Nus88,GG04}.
A similar problem has also been studied in the infinite dimensional setting, when dealing
with Hamilton-Jacobi type PDE: equations arising from deterministic optimal control problems were dealt with in~\cite{Ari97};
the understanding of ergodicity conditions for PDE arising from differential games is a difficult question of current interest, see~\cite{Car10}.

\paragraph{Main result}
In the present paper, we address the question of the solvability of the ergodic equation for zero-sum game problems with {\em finite} state space, but with arbitrary action spaces. 
Our main result, Theorem~\ref{thm:Main}, provides a general sufficient condition for the solvability of the ergodic equation, involving an accessibility condition in directed hypergraphs.
These hypergraphs only depend on the ``behavior at infinity'' of the dynamic programming operator of the game.
Intuitively, a hyperarc encodes the possibility for one of the two players to force the access in one step, with positive probability, to a given set of states, independently of the action of the other player. 
When the action spaces are finite, Theorem~\ref{thm:Main} has a simple game theoretic interpretation:
it shows that the ergodic equation is solvable for all perturbations of the payment function if, and only if, we cannot find two disjoint subsets of states $I,J$ such that Player \MAX, starting from any
state of $J$, has a possibility to force the next state to avoid being in $I$, with positive probability, whereas Player \MIN, starting from any state of $I$, has a possibility to force the next state to avoid being in $J$, still with positive probability.

\paragraph{Related results}
The present results refine or generalize several earlier ergodicity conditions.
In particular, the present hypergraph condition (Theorem~\ref{thm:Main}) refines the ``generalized
Perron-Frobenius theorem'' of Gaubert and Gunawardena~\cite{GG04}, which gave a sufficient condition for the solvability of the ergodic equation in terms of graph accessibility (hypergraphs lead to tighter
conditions).
It also extends results of Cavazos-Cadena and Hern{\'a}ndez-Hern{\'a}ndez~\cite{CCHH10}, which apply to a subclass of operators arising from risk sensitive control problems.
It finally refines a recent result of the authors~\cite{AGH15}, in which reachability conditions involving more special hypergraphs were introduced for games with {\em bounded payments}.
The present Theorem~\ref{thm:Main} is more general, since arbitrary payments are allowed, and when the payments are bounded, the main result of~\cite{AGH15} can be recovered as a special case of the present one, as the hypergraph conditions of Theorem~\ref{thm:Main} reduce to the ones in~\cite{AGH15}.
The case of {\em unbounded payments} is very common when the state space is not compact.
Nevertheless, examples with unbounded payments and finite state space arise typically in the modeling of safety issues, where one player controls the transition probability to an undesirable state with a cost blowing up as this probability vanishes.

The proof of our main result further develops an idea of~\cite{GG04}:
checking the boundedness of the so called ``slice spaces'' with respect to Hilbert's seminorm.
In this respect, we solve one question left open in~\cite{GG04} (characterizing combinatorially the boundedness of all slice spaces).

The paper is organized as follows.
In Section~\ref{sec-prelim}, we recall some basic definitions and results concerning zero-sum stochastic games.
In Section~\ref{sec-main}, we state and establish the main result (ergodicity condition, Theorem~\ref{thm:Main}).
This section includes a discussion of algorithmic and computational complexity aspects.
In Section~\ref{sec-compar}, we compare Theorem~\ref{thm:Main} and its corollary with earlier
results, and provide an illustrative example.

\section{Preliminaries on zero-sum stochastic games}
\label{sec-prelim}

\subsection{The mean payoff}
In this paper, we consider zero-sum stochastic games with {\em perfect information} and finite state space.
We will denote by
$\state := \{1,\dots,n\}$ the state space, and by $\Delta(\state)$ the set of probability measures on $\state$.
For each state $i \in \state$, we denote by $A_i$ (resp.\ $B_i$) the action space of Player \MIN\ (resp.\ \MAX) in state $i$;
$r_i^{a b} \in \R$ the transition payment, paid by Player \MIN\ to Player \MAX\ when the current state is $i$ and the last actions of the players are $a \in A_i$ and $b \in B_i$, respectively;
$P_i^{a b} \in \Delta(\state)$ the transition probability from state $i$ to the next state, assuming that the actions $a \in A_i$ and $b \in B_i$ were just chosen.

Starting from a given initial state $i_0$, a zero-sum game with perfect information is played in stages as follows:
at step $\ell$, if the current state is $i_\ell$, Player \MIN\ chooses an action $a_{i_\ell} \in A_{i_\ell}$, and Player \MAX\ subsequently chooses an action $b_{i_\ell} \in B_{i_\ell}$.
Then, Player \MIN\ pays $r_{i_\ell}^{a_\ell b_\ell}$ to Player \MAX\ and the next state is chosen according to the probability law $P_{i_\ell}^{a_\ell b_\ell}$.
By perfect information, we mean that both players are informed of the current state and of the previous actions of the other player.
In particular, we assumed that Player \MIN\ selects each action before Player \MAX, who is therefore informed of this action.

A strategy $\sigma$ of Player \MIN\ (resp.\ $\tau$ of Player \MAX) is a selection rule determining the current action of Player \MIN\ (resp.\ \MAX) as a function of the current and past information available to him (resp.\ her).
Given strategies $\sigma$ and $\tau$ of the two players, the {\em payoff} of the {\em game in horizon $k$} is the following additive function of the transition payments:
\begin{equation*}
  \label{eq:FiniteHorizonPayoff}
  J_i^k(\sigma,\tau) = \mathbb{E}_{i,\sigma,\tau} \bigg[ \sum_{\ell=0}^{k-1} r_{i_\ell}^{a_\ell b_\ell} \bigg] \enspace ,
\end{equation*}
where $\mathbb{E}_{i,\sigma,\tau}$ denotes the expectation for the probability law of the process $(i_\ell,a_\ell,b_\ell)_{\ell \geq 0}$, determined by the transition probabilities, the strategies $\sigma$ and $\tau$ of the players, and the initial state $i_0 = i$.
Player \MIN\ intends to choose a strategy minimizing the payoff $J_i^k$, whereas Player \MAX\ intends to maximize the same payoff.
The value of the $k$-stage game (played in horizon $k$) starting at state $i$ is thus defined as
\[
  v_i^k = \inf_\sigma \sup_\tau J_i^k(\sigma,\tau) \enspace ,
\]
the infimum and the supremum being taken over the set of strategies of the players \MIN\ and \MAX, respectively.
Note that here the infimum and the supremum commute, since the information is perfect.


The study of the value vector $v^k = (v_i^k)_{i \in \state}$ involves the  {\em dynamic programming operator}, or {\em Shapley operator}.
The latter is a map $T: \R^n \to \R^n$ whose $i$th coordinate is given by
\begin{equation}
  \label{eq:ShapleyOperator}
  T_i(x) = \inf_{a \in A_i} \sup_{b \in B_i} \left( r_i^{a b} + P_i^{a b}x \right) \enspace , \quad x \in \R^n \enspace.
\end{equation}
Note that an element $P \in \Delta(\state)$ is seen as a row vector $P=(P_j)_{j \in \state}$ of $\R^n$, so that $P x$ means $\sum_{j \in \state} P_j x_j$.
The Shapley operator allows one to determine recursively the value vector:
$v^k = T(v^{k-1})$ and $v^0 = 0$.
See~\cite{FV97,NS03,MSS14} for more background.

Here, we are interested in the asymptotic behavior of the sequence of mean values per time unit $(v^k/k)$ as $k\to \infty$.
When the latter ratio converges, the limit, denoted by $\chi(T)$, is called the {\em mean payoff vector}.
Note that, according to the dynamic programming principle, the mean payoff vector is given by
\[
  \chi(T) = \lim_{k \to \infty} \frac{T^k(0)}{k} \enspace ,
\]
where $T^k := T \circ \dots \circ T$ denotes the $k$th iterate of $T$.

The limit does exist for important classes of games.
This question was studied in now classical papers of Bewley and Kohlberg~\cite{BK76} and Mertens and Neyman~\cite{MN81}.
See also Rosenberg and Sorin~\cite{RS01a}, Sorin~\cite{Sor04}, Renault~\cite{Ren11}, and Bolte, Gaubert and Vigeral~\cite{BGV13} for more recent developments.
The existence of $\chi(T)$ is not guaranteed in general: a recent result of Vigeral~\cite{Vig13} shows that the limit may not exist even if the action spaces are compact and the transition payments and probabilities are continuous functions of the actions.

\subsection{Shapley operators}
Given a zero-sum stochastic game with finite state space, its Shapley operator $T$, represented in~\eqref{eq:ShapleyOperator}, satisfies the following two properties:
\begin{itemize}
  \item additive homogeneity: \enspace $T(x+\lambda \unit) = T(x) + \lambda \unit, \enspace \lambda \in \R$;
  \item monotonicity: \enspace $x \leq y \implies T(x) \leq T(y)$;
\end{itemize}
where $\R^n$ is endowed with its usual partial order, and $\unit$ is the {\em unit vector} of $\R^n$, i.e., the vector with all its entries equal to $1$.
The importance of these axioms in stochastic control and game theory has been known for a long time~\cite{CT80}.

Kolokoltsov showed in~\cite{Kol92} that, conversely, every operator from $\R^n$ to $\R^n$ that satisfies the latter two properties can be written as a Shapley operator~\eqref{eq:ShapleyOperator}.
Rubinov and Singer~\cite{RS01b} showed that the game may even be required to be deterministic.
This motivates the following definition.
\begin{definition}
  We call {\em Shapley operator} over $\R^n$ any operator from $\R^n$ to $\R^n$ that is both monotone and additively homogeneous.
\end{definition}

\subsection{The ergodic eigenproblem}
Given a Shapley operator $T$ over $\R^n$, a basic question is to understand when the mean payoff vector is independent of the initial state.
In other words, we ask whether $[\chi(T)]_i = [\chi(T)]_j$ holds for all $i, j \in \state$, i.e., $\chi(T)$ is proportional to the {\em unit vector} $\unit$ of $\R^n$.

This question is related to the nonlinear spectral problem consisting in finding an eigenpair $(\lambda,u) \in \R \times \R^n$ solution of the {\em ergodic equation}
\begin{equation}
  \label{eq:ErgodicEquation}
  T(u) = \lambda \unit + u \enspace .
\end{equation}
Indeed, when the ergodic eigenproblem~\eqref{eq:ErgodicEquation} has a solution, it is easy to see that $\chi(T) = \lambda \unit$ (this follows for instance from Proposition~\ref{prop-elem} below). 
Hence, the scalar $\lambda$, called the {\em eigenvalue} of $T$, yields the mean payoff per time unit for any initial state.
Moreover, the vector $u$, called {\em bias vector} or {\em eigenvector}, gives optimal, or $\varepsilon$-optimal, stationary strategies by identifying in~\eqref{eq:ShapleyOperator} the actions for which the infimum (resp.\ the supremum) is attained, or approached with an arbitrarily small precision.

\subsection{Boundedness of slice spaces in Hilbert's seminorm}
Our approach of the ergodic eigenproblem, following~\cite{GG04}, relies on the study of the {\em sub-eigenspaces} $\Scal_\alpha(T)$, {\em super-eigenspaces} $\Scal^\beta(T)$ and {\em slice spaces} $\Scal_\alpha^\beta(T)$, for $\alpha, \beta \in \R$, which are defined respectively by:
\begin{equation}
  \begin{aligned}
    \Scal_\alpha(T) &= \{ x \in \R^n \mid T(x) \geq \alpha \unit + x \} \enspace , \\
    \Scal^\beta(T) &= \{ x \in \R^n \mid T(x) \leq \beta \unit + x \} \enspace , \\
    \Scal_\alpha^\beta(T) &= \Scal_\alpha(T) \cap \Scal^\beta(T) \enspace .
  \end{aligned}
  \label{e-def-slices}
\end{equation}
These subsets are of intrinsic interest.
Indeed, the following standard proposition shows that by checking the nonemptyness of any of these sets, we can bound the mean payoff from above, from below, or in both directions.
\begin{proposition}[{\cite[Prop.~7]{Sor04}}]
  \label{prop-elem}
  Assume that $\Scal_\alpha(T) \neq \emptyset$. Then, $\liminf_{k \to \infty} T^k(0) / k \geq \alpha \unit$.
  Similarly, if $\Scal^\beta(T) \neq \emptyset$, then, $\limsup_{k \to \infty} T^k(0) / k \leq \beta \unit$.
\end{proposition}
Therefore, an element $u$ such that $T(u) \geq \alpha \unit + u$ may be thought of as a {\em certificate} that $\liminf_{k \to \infty} T^k(0) / k \geq \alpha \unit$.
A dual interpretation holds when $T(u) \leq \beta \unit + u$.
Also, if the ergodic equation~\eqref{eq:ErgodicEquation} has a solution, we deduce that $\Scal_\lambda^\lambda \neq \emptyset$, and so, it follows from Proposition~\ref{prop-elem} that $\chi(T)=\lambda \unit$. 

Let us now introduce the {\em Hilbert's seminorm} on $\R^n$, which is given by
$ \Hnorm{x} = \max_{i \in \state} x_i - \min_{i \in \state} x_i$.
It is a useful tool in nonlinear Perron-Frobenius theory~\cite{Nus88,GG04}, where it is sometimes called {\em Hopf's oscillation}~\cite{Hop63,Bus73}.
It also appears in linear consensus theory under the name of {\em diameter} or {\em Tsitsiklis' Lyapunov function}~\cite{TBA86}.
In this context, the set $\R \unit$ of scalar multiples of the unit vector is known as {\em consensus states}, and the Hilbert's seminorm, which vanishes on $\R \unit$, is a measure of the distance to consensus states. 
See~\cite{SSR10,GQ13} for further developments concerning Hilbert's seminorm and the related Hilbert's projective metric in relation with consensus theory. 

It is a standard result (see e.g.~\cite{GG04}) that any Shapley operator $T: \R^n \to \R^n$ is nonexpansive with respect to Hilbert's seminorm, meaning that, for every $x, y \in \R^n$,
$\Hnorm{T(x)-T(y)} \leq \Hnorm{x-y}$.
Thus, as a special case of Theorem~4.1 in~\cite{Nus88}, we get that the ergodic eigenproblem~\eqref{eq:ErgodicEquation} has a solution if, and only if, one orbit of $T$, $\{ T^k(x) \mid k \in \N \}$ with $x \in \R^n$, is bounded in Hilbert's seminorm.
Hence, the solvability of~\eqref{eq:ErgodicEquation} boils down to finding a nonempty invariant set under $T$ that is bounded in Hilbert's seminorm. 
It is readily seen that all the subsets in~\eqref{e-def-slices} are invariant under $T$.
Then, we have the following result.
\begin{theorem}[{\cite[Th.~4.1]{Nus88}, \cite{GG04}}]
  \label{thm:BoundedSliceSpace}
  Let $T: \R^n \to \R^n$ be a Shapley operator.
  Assume that there is a choice of $\alpha, \beta \in \R$ such that the slice space $\Scal_\alpha^\beta(T)$ is nonempty and bounded in Hilbert's seminorm.
  Then, there exist $\lambda \in \R$ and $u \in \R^n$ such that $T(u) = \lambda \unit + u$.
\end{theorem}
Note that the slice space $\Scal_\alpha^\beta(T)$ is nonempty as soon as $\alpha \leq \min_{i \in \state} T_i(0)$ and $\beta \geq \max_{i \in \state} T_i(0)$, since in this case $0 \in \Scal_\alpha^\beta$.
Hence, the main difficulty when applying Theorem~\ref{thm:BoundedSliceSpace} is to check that the boundedness condition does hold for some specific $\alpha, \beta$.
We shall see in the next section that the boundedness of the slices spaces {\em for all} values of $\alpha, \beta$ has a combinatorial characterization.
Also, we easily deduce the following result.
\begin{corollary}
  \label{coro:BoundedSliceSpaces}
  Let $T: \R^n \to \R^n$ be a Shapley operator.
  If all the slice spaces $\Scal_\alpha^\beta(T)$ are bounded in Hilbert's seminorm, then, for every vector $g \in \R^n$, there is a scalar $\lambda \in \R$ and a vector $u \in \R^n$ such that $g+T(u) = \lambda \unit + u$.
\end{corollary}

\section{Hypergraph conditions for the solvability of the ergodic eigenproblem}
\label{sec-main}

\subsection{Hypergraphs associated with a Shapley operator}
\label{HypergraphShapleyOperator}
A {\em directed hypergraph} is a pair $(N,D)$, where $N$ is a set of {\em nodes} and $D$ is a set of (directed) {\em hyperarcs}.
A hyperarc $d$ is an ordered pair $(\tail(d),\head(d))$ of disjoint nonempty subsets of nodes;
$\tail(d)$ is the {\em tail} of $d$ and $\head(d)$ is its {\em head}.
We shall often write for brevity $\tail$ and $\head$ instead of $\tail(d)$ and $\head(d)$, respectively.
When $\tail$ and $\head$ are both of cardinality one, the hyperarc is said to be an arc, and when every hyperarc is an arc, the directed hypergraph becomes a directed graph.
In the following, the term {\em hypergraph} will always refer to a directed hypergraph.
For background on hypergraphs, we refer the reader to~\cite{All14} and the references therein, and in particular to~\cite{GLNP93} for reachability problems.

To express the ``behavior at infinity'' of a Shapley operator $T: \R^n \to \R^n$, we introduce a pair of hypergraphs associated with $T$, and denoted by $(\Hcal^+,\Hcal^-)$.
Before giving their construction, let us fix some notation.
We still denote by $\state$ the set $\{1,\dots,n\}$, and if $I$ is a subset of $\state$, we denote by $\unit_I$ the indicator vector of $I$, i.e., the vector of $\R^n$ such that $[\unit_I]_i = 1$ if $i \in I$ and $[\unit_I]_i = 0$ if $i \in \cpt{I}$.
If $I=\{i\}$, we will write $\unit_i$ instead of $\unit_{\{i\}}$.

We now explain the construction of the hypergraphs $\Hcal^+$ and $\Hcal^-$ associated with $T$:
\begin{itemize}
  \item the set of nodes for both $\Hcal^+$ and $\Hcal^-$ is $\state$;
  \item the hyperarcs of $\Hcal^+$ are the pairs $(J,\{i\})$ such that $\lim_{\alpha \to +\infty} T_i(\alpha \unit_J) = +\infty$;
  \item the hyperarcs of $\Hcal^-$ are the pairs $(J,\{i\})$ such that $\lim_{\alpha \to -\infty} T_i(\alpha \unit_J) = -\infty$.
\end{itemize}
These hypergraphs have an intuitive game theoretic interpretation when the action spaces $A_i$ and $B_i$ are finite.
Then, it can be readily checked that there is a hyperarc $(J,\{i\})$ in $\Hcal^-$ if, and only if, Player \MIN\ can choose an action $a \in A_i$ in such a way that for all actions $b \in B_i$ of Player \MAX, $[P^{ab}_i]_j>0$ holds for at least one state $j \in J$.
In other words, such a hyperarc encodes the fact that starting from state $i$, Player \MIN\ can force the next state to belong to $J$, with positive probability.
Similarly, there is a hyperarc $(J,\{i\})$ in $\Hcal^+$ if, and only if, for every action $a \in A_i$, Player \MAX\ can choose an action $b \in A_i$ (depending on $a$) such that $[P_i^{ab}]_j>0$ holds for at least one state $j \in J$.
Therefore, starting from state $i$, Player \MAX\ can force the next state to belong to $J$, with positive probability.

\subsection{Reachability conditions for bounded slice spaces}
We first recall the definition of reachability in hypergraphs.
Let $\Gcal=(N,D)$ be a hypergraph.
A {\em hyperpath} of length $p$ from a set of nodes $I \subset N$ to a node $j \in N$ is a sequence of $p$ hyperarcs $(\tail_1,\head_1),\dots,(\tail_p,\head_p)$, such that $\tail_i \subset \bigcup_{k=0}^{i-1} \head_k$ for all $i=1,\dots,p+1$ with the convention $\head_0=I$ and $\tail_{p+1}=\{j\}$.
We say that a node $j \in N$ is {\em reachable} from a subset $I$ of $N$ if there exists a hyperpath from $I$ to $j$.
Alternatively, the relation of reachability can be defined in a recursive way:
$j$ is reachable from $I$ if either $j \in I$ or there exists a hyperarc $(\tail,\head)$ such that $j \in \head$ and every node of $\tail$ is reachable from the set $I$.
A subset $J$ of $N$ is said to be {\em reachable} from a subset $I$ of $N$ if every node of $J$ is reachable from $I$.
We denote by $\reach(I,\Gcal)$ the set of reachable nodes from $I$ in $\Gcal$.
A subset $I$ of $N$ is {\em invariant} in the hypergraph $\Gcal$ if it contains every node that is reachable from itself, i.e., $\reach(I,\Gcal) \subset I$ (hence $\reach(I,\Gcal) = I$, since the other inclusion always holds).
One readily checks that, for $J \subset N$, $\reach(J,\Gcal)$ is the smallest invariant set in $\Gcal$ containing $J$.


We also need the following definition.
\begin{definition}
  \label{def:ConjugateSubsets}
  A pair $(I,J)$ of subsets of $\state$ is said to be a pair of {\em conjugate subsets} with respect to the hypergraphs $(\Hcal^+,\Hcal^-)$ if $I \cap J = \emptyset$, $\reach(J,\Hcal^+) = \cpt{I}$ and $\reach(I,\Hcal^-) = \cpt{J}$.
  Such a pair is said to be {\em trivial} if one of the sets $I,J$ is empty.
\end{definition}

The following is our main result.
\begin{theorem}
  \label{thm:Main}
  Let $T: \R^n \to \R^n$ be a Shapley operator.
  All the slice spaces $\Scal_\alpha^\beta(T)$ are bounded in Hilbert's seminorm if, and only if, there are only trivial pairs of conjugate subsets with respect to the hypergraphs $(\Hcal^+, \Hcal^-)$. 
\end{theorem}

The proof will be given in the extended version of the present conference article.
Let us just mention that the ``if'' part is proved by combinatorial means, and the ``only if'' part by a topological degree theory argument. 

Combining this theorem and Corollary~\ref{coro:BoundedSliceSpaces}, we get the following existence result for solutions of the ergodic eigenproblem~\eqref{eq:ErgodicEquation}.
\begin{corollary}
  \label{coro:ErgodicityCondition}
  Let $T: \R^n \to \R^n$ be a Shapley operator.
  Assume that there are only trivial pairs of conjugate subsets with respect to the hypergraphs $(\Hcal^+, \Hcal^-)$.
  Then, for every vector $g \in \R^n$, there is a scalar $\lambda \in \R$ and a vector $u \in \R^n$ such that $g+T(u) = \lambda \unit + u$.
\end{corollary}

It has already been mention in the introduction (see also Subsection~\ref{subsec-recent}), that the hypergraph conditions of Corollary~\ref{coro:ErgodicityCondition} reduce to the ones in~\cite{AGH15} when the payments are bounded.
Thus, the main result of~\cite{AGH15} implies that the converse property in Corollary~\ref{coro:ErgodicityCondition} does hold if, for every state $i \in \state$, the payment function $(a,b) \mapsto r_{i}^{ab}$ is bounded.
In a forthcoming work, we show that the converse property is in fact true with arbitrary payments.

\subsection{Algorithmic aspects}
Given a Shapley operator $T$ over $\R^n$, the basic issue under consideration is to check whether the ergodic eigenproblem~\eqref{eq:ErgodicEquation} is solvable for all operators $g+T$ with $g \in \R^n$.

Corollary~\ref{coro:ErgodicityCondition} provides a combinatorial 
condition for this property to hold.
This condition can be effectively checked as soon as the limits $\lim_{\alpha \to \pm \infty} T_i(\alpha \unit_J)$ arising in the definition of the hyperarcs of $\Hcal^\pm$ can be computed. 
This is generally the case in practice, as the Shapley operator is typically given explicitly in terms of actions, payments, and transition probabilities. 
To give an elementary example, we already noted at the end of Subsection~\ref{HypergraphShapleyOperator} that if all the action spaces $A_i, B_i$ in~\eqref{eq:ShapleyOperator} are {\em finite}, then, checking whether there is a hyperarc $(J,\{i\})$ in $\Hcal^\pm$ reduces to elementary tests.
Indeed, the discussion in that subsection implies that $(J,\{ i\})$ is a hyperarc of $\Hcal^-$ (resp.\ $\Hcal^+$) if, and only if
\[
  \max_{a \in A_i} \min_{b \in B_i} P^{ab}_{i} \unit_J >0 \quad (\text{resp. } \min_{a \in A_i} \max_{b \in B_i} P^{ab}_{i}\unit_J >0) \enspace .
\]
Therefore, checking whether the hyperarc $(J,i)$ belongs to $\Hcal^\pm$ can be done in $O(n \, |A_i| \, |B_i|)$ elementary operations, where $|\cdot|$ denotes the cardinality of a set. 
However, the limits $\lim_{\alpha \to \pm \infty} T_i(\alpha \unit_J)$ are computable in more general situations (see Section~\ref{ex-shapley2}), and so the interest of Corollary~\ref{coro:ErgodicityCondition} is to lead to an effective 
condition in such cases.

Let us now recall how reachability conditions can be checked algorithmically once the hypergraphs have been computed.
Indeed, the conditions of Corollary~\ref{coro:ErgodicityCondition} boil down to check that, for every $I \subset \state$, either $\reach(\cpt{I},\Hcal^+) \neq \cpt{I}$ or $\reach(I,\Hcal^-) = \state$.

It is known that in a hypergraph $\Gcal=(N,D)$, the set of reachable nodes from a subset $I$ can be computed in $O(\size(\Gcal))$ time, where $\size(\Gcal)$ denotes the {\em size} of $\Gcal$, and is equal to $|N| + \sum_{d \in D} |\tail(d)| + |\head(d)|$, see~\cite{GLNP93}.
Hence, the conditions of Corollary~\ref{coro:ErgodicityCondition} can be checked in time $O(2^n \size(\Hcal^\pm))$, counting for one time unit a call to an oracle computing a limit $\lim_{\alpha\to\pm\infty}T_i(\alpha \unit_J)$. 
However, $\size(\Hcal^\pm)$ may be exponential in the number of states, since the number of hyperarcs is only bounded by $n 2^n$. 
The bound $O(2^n \size(\Hcal^\pm))$ may appear to be coarse at first sight, but it cannot be reduced to a polynomial bound unless $\text{P} = \text{coNP}$. 
Indeed, the restricted version of this problem concerning deterministic operators with finite action spaces reduces to the nonexistence of a nontrivial fixed point of a monotone Boolean function, a problem which has been shown to be coNP-hard by Yang and Zhao~\cite{YZ04}.

We refer the reader to~\cite{AGH15} for more information on complexity issues, and refinements in the case of finite actions spaces (in this case, one arrives at better complexity bounds, still exponential in $n$).
We also point out that there are simpler sufficient conditions, involving directed graphs instead of a hypergraph~\cite{GG04}.
Such conditions can be checked using only a polynomial number of elementary operations (still counting
for one unit the call to the oracle computing a limit), but they are less accurate.

\section{Comparison with earlier results and illustrative example}
\label{sec-compar}

We first show that several earlier results follow as special cases of Theorem~\ref{thm:Main}.
Then, we give 
a simple example in which Theorem~\ref{thm:Main} allows one to show the solvability of the ergodic problem, whereas earlier conditions fail. 

\subsection{The ergodic eigenproblem and the recession operator}
To address the problem of the solvability of the ergodic equation~\eqref{eq:ErgodicEquation}, Gaubert and Gunawardena introduced in~\cite{GG04} the {\em recession operator} associated with a Shapley operator $T$ over $\R^n$.
This operator is defined, for $x \in \R^n$, by
$\hat{T}(x) = \lim_{\alpha \to +\infty} T(\alpha x) / \alpha$.
Its existence is not always guaranteed, but it does exist for any Shapley operator represented as in~\eqref{eq:ShapleyOperator} with bounded payments.
Moreover, when it exists, it inherits from $T$ the additive homogeneity and monotonicity properties. 
Furthermore, it is positively homogeneous, meaning that $\hat{T}(\alpha x) = \alpha \hat{T}(x)$ for every $\alpha \geq 0$.
As a consequence, any vector proportional to the unit vector of $\R^n$ is a fixed point of $\hat T$.
We shall call such fixed points {\em trivial} fixed points.

Theorem~13 of~\cite{GG04} shows that if a Shapley operator $T$ has a recession function $\hat{T}$ whose only fixed points are trivial, then all the slice spaces of $T$ are bounded in Hilbert's seminorm.
An example given in the same paper shows that this sufficient condition is not necessary.
Therefore, the present Theorem~\ref{thm:Main} solves the issue raised there of giving a necessary and sufficient condition.

\subsection{Weakly convex maps and risk sensitive control}

Cavazos-Cadena and Hern{\'a}ndez-Hern{\'a}ndez gave in~\cite{CCHH10}, under a weak convexity property, a necessary and sufficient condition, in terms of graph, under which the ergodic eigenproblem~\eqref{eq:ErgodicEquation} is solvable for all the Shapley operators $g+T$ with $g \in \R^n$.
They also showed that this condition is equivalent to $\hat{T}$ having only trivial fixed points.

The weak convexity property was motivated by applications to risk sensitive control.
A typical example of weakly convex operator is the following Shapley operator $T$:
\begin{equation}
  T_i(x) = \log \Big( \sum_{j \in \state} M_{ij} \exp(x_j) \Big) \enspace, \quad i \in \state \enspace ,
  \label{e-risk}
\end{equation}
where $M=(M_{ij})$ is a nonnegative matrix without zero row.
A supremum of weakly convex Shapley operators is weakly convex, hence, one can construct further examples of weakly convex operators by taking suprema of operators of the form~\eqref{e-risk}.

It can be easily checked that under the weakly convex hypothesis, our hypergraph conditions (Theorem~\ref{thm:Main}) reduce to the graph conditions of~\cite{CCHH10}.
However, this weak convexity assumption is typically not satisfied for two-player Shapley operators, or for dynamic programming operators of general stochastic control problems.

\subsection{Hypergraphs associated with a recession operator}
\label{subsec-recent}
A recent result of the authors~\cite{AGH15} concerns the special case in which the payment functions $(a,b) \mapsto r_i^{ab}$ are {\em bounded} for all $i \in \state$.
Then, the recession map of the Shapley operator $T$ exists and each of its coordinate functions is given by 
\begin{equation*}
  \label{eq:PaymentFreeOperator}
  \hat T_i(x) = \min_{a \in A_i} \max_{b \in B_i} P_i^{a b}x \enspace , \quad i \in \state \enspace .  
\end{equation*}
Theorem~3.1 of~\cite{AGH15} shows that the ergodic equation is solvable for all operators $g+T$ with $g \in \R^n$ if, and only if, the recession operator $\hat T$ has only trivial fixed points.
It is also shown there that the latter condition can be characterized in terms of reachability in directed hypergraphs, depending only on the supports of the probability laws $P_i^{ab}$. 
We leave it as an exercise to check that the hypergraph conditions of Theorem~\ref{thm:Main} reduces to the one of~\cite{AGH15} in this special case.
It follows that the sufficient condition of Corollary~\ref{coro:ErgodicityCondition} is necessary when the payment functions are bounded.


\subsection{Example}
\label{ex-shapley2}

Here, we illustrate the situation in which the conditions for the solvability of the eigenproblem~\eqref{eq:ErgodicEquation} given by Corollary~\ref{coro:ErgodicityCondition} are satisfied, whereas the conditions provided by~\cite{GG04,CCHH10,AGH15} do not apply.
Consider the parametric family of Shapley operators $g+T$, $g \in \R^3$, with $T: \R^3 \to \R^3$ given by
\[
  T(x) =
  \begin{pmatrix}
    \sup_{0 < p \leq 1} \big( \log p + p ( x_2 \wedge x_3 ) + (1-p) x_1 \big) \\
    \inf_{0 < p \leq 1} \big( -\log p + p x_3 + (1-p) x_1 \big) \\
    x_3
  \end{pmatrix}
\]
where $\wedge$ stands for $\min$.
The operator $g+T$ corresponds to a game with three states:
Player \MAX\ partially control state $1$, Player \MIN\ controls state $2$ and state $3$ is an absorbing state.
In state $1$, Player \MAX\ first receives $g_1$ from Player \MIN, then chooses an action $p \in (0,1]$, and receives in addition $\log p$.
Then, with probability $1-p$, the next state remains $1$, and with probability $p$, it is chosen by Player \MIN\ between state $2$ and state $3$.
Thus, maximizing the one day payoff would lead to select $p=1$, but this leads to leave state $1$ with probability one.
A dual interpretation applies to Player \MIN in state $2$.

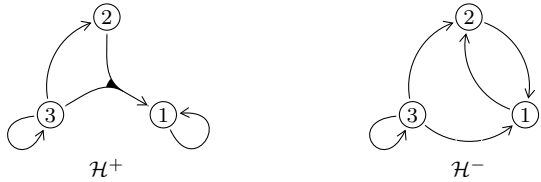
\begin{figure}[!h]
  \begin{minipage}[b]{0.49\linewidth}
    \centering \footnotesize
    \begin{tikzpicture}
      [on grid,auto,bend angle=40,
      state/.style={circle,draw,inner sep=0pt,minimum size=0.35cm},
      phantom/.style={draw,inner sep=0pt,minimum size=0.01ex},
      head/.style={<-,shorten <=1pt,shorten >=1pt,>=angle 60},
      tail/.style={->,shorten <=1pt,shorten >=1pt,>=angle 60},
      loopstyle/.style={->,shorten <=1pt,shorten >=1pt,>=angle 60,min distance=0.75cm,looseness=10},
      hyperedge/.style={>=angle 60}]
      \node[state] (s1) at (0,0) {$1$}
        edge [loop,loopstyle,in=5,out=-65] (s1);
      \node[state] (s2) at (-0.75,1.3) {$2$};
      \node[state] (s3) at (-1.5,0) {$3$}
        edge [tail,bend left] (s2)
        edge [loop,loopstyle,in=-115,out=185] (s3);
      \node[phantom] (s11) at (-0.2,0.12) {};
      \node[phantom] (s22) at (-0.75,1.07) {};
      \node[phantom] (s33) at (-1.3,0.12) {};
      \hyperedgewithangles[0.3][$(hyper@tail)!0.6!(hyper@head)$]{s22/-90,s33/30}{-30}{s11/150};
      \node[state,color=black!0] (s5) at (-0.75,-0.3) [label=below:{$\Hcal^+$}] {};
    \end{tikzpicture}
  \end{minipage}
  \hfill
  \begin{minipage}[b]{0.49\linewidth}
    \centering \footnotesize
    \begin{tikzpicture}
      [on grid,bend angle=40,auto,
      state/.style={circle,draw,inner sep=0pt,minimum size=0.35cm},
      head/.style={<-,shorten <=1pt,shorten >=1pt,>=angle 60},
      tail/.style={->,shorten <=1pt,shorten >=1pt,>=angle 60},
      loopstyle/.style={->,shorten <=1pt,shorten >=1pt,>=angle 60,min distance=0.75cm,looseness=10}]
      \node[state] (s1) at (0,0) {$1$};
      \node[state] (s2) at (-0.75,1.3) {$2$}
        edge [tail,bend left] (s1)
        edge [head, bend right] (s1);
      \node[state] (s3) at (-1.5,0) {$3$}
        edge [loop,loopstyle,in=-115,out=185] (s3)
        edge [tail,bend right] (s1)
        edge [tail,bend left] (s2);
      \node[state,color=black!0] (s5) at (-0.75,-0.3) [label=below:{$\Hcal^-$}] {};
    \end{tikzpicture}
  \end{minipage}
  \caption{The hypergraphs associated with $T$}
  \label{fig:Hpm}
\end{figure}

Figure~\ref{fig:Hpm} shows a concise representation of the hypergraphs $\Hcal^+$ and $\Hcal^-$ associated with $T$, in which only the hyperarcs with minimal tail (with respect to the inclusion partial order) have been represented.
To construct these hypergraphs, it is convenient to notice that $T_1(x) = h((x_2 \wedge x_3) - x_1) + x_1$ and $T_2(x) = -h(x_1-x_3)+x_1$, where $h$ is the real function defined by $h(z) = \sup_{0 < p \leq 1} (\log p + p z)$.
Also note that $h$ satisfies $h(z) = -1-\log(-z)$ for $z \leq -1$, and $h(z)=z$ for $z \geq -1$.
Thus, for instance, there is no arc from $\{2\}$ to $\{1\}$ in $\Hcal^+$ since $T_1(\alpha \unit_2) = 0$ for all $\alpha \geq 0$.
However, there is a hyperarc from $\{2,3\}$ to $\{1\}$, since $T_1(\alpha \unit_{\{2,3\}}) = \alpha$ for all $\alpha \geq 0$, which yields $\lim_{\alpha \to +\infty} T_1(\alpha \unit_{\{2,3\}}) = +\infty$.
Alternatively, following the discussion at the end of Subsection~\ref{HypergraphShapleyOperator}, one may checked that, starting from state $1$, Player \MAX\ can force the next state to belong to $\{2,3\}$ with positive probability, but he cannot force it to be $2$.

Then, one may check that there are no nonempty disjoint subsets of states $I,J$ such that $\reach(J,\Hcal^+) = \{1,2,3\} \setminus I$ and $\reach(I,\Hcal^-) = \{1,2,3\} \setminus J$.
Hence, there do not exist nontrivial conjugate subsets with respect to $(\Hcal^+,\Hcal^-)$.
By application of Theorem~\ref{thm:Main}, we deduce that all the slice spaces of $T$ are bounded in Hilbert's seminorm.
It follows that the ergodic equation~\eqref{eq:ErgodicEquation} is solvable for all operators $g+T$, with $g \in \R^3$.

This conclusion cannot be obtained from the previous theorems in~\cite{GG04,CCHH10,AGH15}.
Indeed, although the recession operator of $T$ exists (it is given by $\hat T_1(x) = x_1 \vee (x_2 \wedge x_3)$, $\hat T_2(x) =  x_1 \wedge x_3$ and $\hat T_3(x) = x_3$), one may check that any vector $[\alpha,0,0]^\top$ with $\alpha \geq 0$ is a fixed point of $\hat T$.
Thus, ergodicity conditions requiring the nonexistence of nontrivial fixed points of the recession operator $\hat T$ cannot be applied.




\bibliographystyle{IEEEtran}
\bibliography{references}

\end{document}